\documentclass[12pt,a4paper]{article}

\usepackage{pdflscape}
\usepackage{amsmath}
\usepackage{amssymb}
\usepackage{latexsym}
\usepackage{srcltx}
\usepackage{graphics}
\usepackage{color}
\usepackage{epsfig}
\usepackage{color}

\newcommand{\re}{{\mathbb R}}

\newcommand{\cA}{{\mathcal{A}}}

\newcommand{\bx}{{\boldsymbol{x}}}
\newcommand{\by}{{\boldsymbol{y}}}
\newcommand{\bz}{{\boldsymbol{z}}}
\newcommand{\be}{{\boldsymbol{e}}}

\newcommand{\ba}{{\boldsymbol{a}}}

\newcommand{\bu}{{\boldsymbol{u}}}

\newcommand{\bh}{{\boldsymbol{h}}}

\newtheorem{theorem}{Theorem}
\newtheorem{prop}{Proposition}

\newtheorem{cor}{Corollary}
\newtheorem{remark}{Remark}
\newtheorem{ex}{Example}
\newtheorem{defi}{Definition}

\date{}

\author{Maxim Makarov 
\thanks{Moscow State University; 
Moscow Center of Fundamental and Applied Mathematics, 
Leninstie Gory 1, Moscow  119991 (Russia); {e-mail: \tt\small
maximka1905@mail.ru}}, 
Vladimir Yu. Protasov 
\thanks{DISIM, University of L'Aquila, Italy; {e-mail: \tt\small
vladimir.protasov@univaq.it}}
}

\title{Autopolar conic bodies and polyhedra
\thanks{The first  author is supported   by the Russian Science Foundation no. 23-71-30001. 
The results of  Sections  2, 4, and 6 are obtained  by the first author; the results of Sections 
2, 3, and 5 are obtained by the second author. 
}}

\begin{document}
\maketitle

\begin{abstract}

An antinorm is a concave analogue of a norm. In contrast to norms, antinorms are 
not defined  on the entire space~$\re^d$ but  
  on a cone~$K\subset \re^d$. They  are applied in the matrix analysis, optimal control,  and dynamical systems. Their level sets are called conic bodies and (in case of piecewise-linear antinorms) conic polyhedra. The basic facts
  and notions  of the ``concave analysis'' of antinorms such as 
  separation theorems, duality, polars, Minkowski functionals, etc., are similar to 
those from the standard convex analysis. There are, however, some 
significant differences. One of them is  
the existence of many self-dual objects. We prove that there are infinitely many families of 
autopolar conic bodies and polyhedra in the cone~$K=\re^d_+$. 
 For~$d=2$, this gives a complete classification of 
 self-dual antinorms, while for~$d\ge 3$, there are counterexamples.

\bigskip

\noindent \textbf{Key words:} {\em antinorm, cone, polyhedron, duality, polar transform, 
autopolar set, self-duality
}
\smallskip

\begin{flushright}
\noindent  \textbf{AMS 2020 subject classification} {\em 52A21, 52B11, 46B10}

\end{flushright}

\end{abstract}
\bigskip

\vspace{1cm}

\begin{center}

\large{\textbf{1. Introduction}}	
\end{center}
\bigskip

\begin{center}

\textbf{1.1. Antinorms and conic polyhedra}	
\end{center}
\bigskip 

The word  ``antinorm'' has been understood  in different senses. 
Following most of the  literature we define it as a ``concave norm'', i.e.,  a nonnegative, positively homogeneous  concave function. It is shown easily that such a function
cannot be defined on the entire~$\re^d$, unless it is an identical zero. 
Therefore, an antinorm is usually restricted to a cone. 
Brief historical remarks  and examples of applications are given below in subsection~1.2. 
The level sets of an antinorm  are convex and unbounded. Such ``antiballs'' 
are referred to as conic bodies; in particular, conic polyhedra are conic bodies bounded by several hyperplanes. Most geometrical properties of those objects are similar to 
what we have  in the classical convex analysis.  
There are, however, some significant differences, one of which, the self-duality phenomenon,   is the subject of this paper.

We consider a cone~$K \subset \re^d$ with an apex at the origin. 
It is  assumed to be closed, solid (with a nonempty interior) and 
pointed (not containing a straight line). As usual, faces 
of a cone are its intersections with planes of support. 
The apex is a face of dimension zero, the face of dimension one is an {\em edge}, a face of dimension~$d-1$ is a {\em facade}.  A cone is polyhedral if 
it is defined by a system of several linear inequalities. 
An {\em interior ray} is a ray going from the origin and passing 
through an interior point. The dual cone is defined in the usual way as~$K^* = \bigl\{\by \in \re^d: \, \inf\limits_{\bx \in K} (\by, \bx) \ge 0\bigr\}$.

\begin{defi}\label{d.10}
An antinorm on a cone~$K$ is a concave homogeneous  nonnegative, somewhere positive,  
 function~$f:K\to \re_+$. 
\end{defi}
Homogeneity means that~$f(\lambda \bx) = \lambda f(\bx), \, \lambda \in \re$. Since 
the cone~$K$ is pointed, it follows that~$\lambda$ must be nonnegative, unless~$\bx =0$, otherwise $\lambda\bx \notin K$. An antinorm may vanish  only on the 
boundary of~$K$. It may not be continuous, but all its points 
of discontinuity are on the boundary. Surprisingly enough, 
there ere examples of antinorms that do not have 
continuous extensions from the interior of~$K$ to the boundary. 
However, such an extension always exists and unique for polyhedral cones, in 
particular, for~$K = \re^d_+$, see~\cite{P22}.

The set~$G=\{\bx \in K: \ f(\bx) \ge R\}$ is a {\em ball of radius~$R$}
of the antinorm~$f$. Since  it is defined by a reverse inequality ($f(\bx) \ge R$
instead of
$f(\bx) \le R$), it should be rather called an ``antiball'', but we use a simplified notation. If~$\bx \in G$, then~$\lambda \bx \in G$ for all~$\lambda \ge 1$, therefore, 
a ball of the antinorm is never bounded. We usually deal with unit balls, when~$R=1$. 
The set~$S = \partial G$ is a {\em sphere} of the antinorm. 
\begin{defi}\label{d.20}
For a given cone~$K$, a {\em conic body}~$G$ is a closed convex proper subset of~$K$ 
such that for every interior ray~$\ell \subset K$, the set~$G\cap \ell$ is a ray. 
\end{defi}
A ball of an antinorm is a conic body. Conversely, every 
conic body is a unit ball of the antinorm~$f_G(\bx) \, = \, 
\sup\, \{\lambda > 0: \ \lambda^{-1}\bx \in G\}$, which is an analogue of the standard Minkowskii functional. 
A finite intersection of conic bodies is a conic body, the same is true for the Miknowski sum. However, this is in general not true 
for the convex hull. This operation is replaced by 
the {\em positive convex hull}
${\rm co}_{+} (X) \, = \, {\rm co} (X) \, + \, K$, where
 ${\rm co}(\cdot)$ denotes the 
standard convex hull. 
For every~$X\subset K$, the closure of~${\rm co}_{+} (X)$ is either a conic 
body or the entire cone~$K$. In particular, the positive convex hull of a 
finite set not containing the origin is a {\em conic polytope}. 
\begin{defi}\label{d.30}
A {\em conic polyhedron}~$G$ is a set~$P = \{\bx \in K: \ 
(\ba_j, \bx) \ge 1, \, j=1, \ldots , n\}$, where~$\ba_j$
are nonzero elements of~$K^*$.  
\end{defi}
Thus, a conic polyhedron is a conic body bounded by several hyperplanes. 
The corresponding antinorm is piecewise-linear: $\, f_P(\bx) \, = \, 
\min_{i} (\ba_i, \bx)$. 
In general, a conic  polyhedron may not be a conic polytope. Nevertheless, if the cone $K$ is polyhedral, 
in particular, if~$K = \re^d_+$, then those two notions coincide. 

For the sake of simplicity, in what follows we assume that~$K=\re^d_+$. Since each antinorm on~$\re^d_+$ has a unique continuous extension, 
we consider only continuous antinorms. 
As we see, the main facts of the classical convex analysis have analogues for 
antinorms and for conic bodies/polyhedra, with replacement of $\sup$ by~$\inf$ 
and of~${\rm co}(\cdot)$ by~${\rm co}_+(\cdot)$. The same can be said of the duality theory
with one exception: 
while in the classical convex analysis there is a unique autopolar set 
(the unit Euclidean ball), the family of autopolar conic bodies 
and polyhedra in~$\re^d_+$ is quite  rich and diverse. We give formal definitions and preliminary facts in the next section and then, in Section~3 we formulate the main results. 
Theorem~\ref{th.20} will allow us to construct infinitely many families of self-dual antinorms and autopolar conic bodies. 
In case~$d=2$ this gives a complete classification 
of autopolar sets (Theorem~\ref{th.30}) while for~$d=3$, there is a counderexample (Theorem~\ref{th.75})
showing that not all such sets are produces by~Theorem~\ref{th.20}. 
A complete classification of self-dual antinorms and of autopolar sets is formulated as an open problem in~Section~6. 
\smallskip 

Throughout the paper we denote vectors by bold letters, their coordinates
and other scalars by usual letters, so,~$\bx = (x_1, \ldots , x_d)$. 
The Euclidean norm is denoted by~$|\bx|$. By {\em projection} we always mean the 
orthogonal projection of~$\re^d$ to a proper affine subspace. 
A preimage of this operator will be referred to as  {\em orthogonal extension}. 
So, the orthogonal extension of a set~$X \subset \re^d$ 
is $\{\bx + \bh: \ \bx \in X,\ \bh\in X^{\perp}\}$, where~$X^{\perp}$
denotes the orthogonal complement to the linear span of~$X$. 

\bigskip

\begin{center}

\textbf{1.3. Applications and related works} 	
\end{center}
\bigskip

To the best of our knowledge, Definition~\ref{d.10} 
 originated with Merikoski~\cite{M90}
 and the first results on antinorms were obtained  in \cite{MS91, MO92}. Matrix antinorms  on the 
 positive semidefinite cone  were 
studied in~\cite{BH11,  BH14}; some works considered the Minkowski antinorm $f(A) = ({\rm det}\, A)^{1/d}$ and 
 the Schatten $q$-antinorms~$f(A) = ({\rm tr}\, A)^{1/q}, \, q \in (-\infty, 1]$, see also~\cite{BH15} for extensions to von Neumann algebras. Independently and later  antinorms were defined 
 in~\cite{P10} for analysing infinite products of random matrices. 

Antinorms have been  applied 
in the study of linear dynamical systems 
of the form \linebreak $\, \dot \bx \, = \, A(t)\bx, \, t\ge 0, $ (continuous time) and 
~$\, \bx_{k+1} = A(k)\bx_k\, , \, k\ge 0$,  (discrete time), where the matrix~$A(\cdot)$
is chosen independently from a given compact set for each value of the argument. If the system is asymptotically stable i.e., 
all its trajectories~$\bx_k$ tend to zero as~$k\to \infty$, then it possesses 
 a convex Lyapunov function ({\em Lyapunov norm}). 
 This norm  is the main tool to prove the stability~\cite{L03, MP89, LA09}. However, for other types of stability, the Lyapunov norm may not exist. 
 This occurs, for instance, for the {\em almost sure stability}, when 
 almost all trajectories tend to zero~\cite{DJ22, H97}. Nevertheless, in this case the Lyapunov function can be constructed  as a {\em concave} homogeneous function i.e., the 
 antinorm~\cite{JP13,  P13, P11}. 
 The same can be said on the {\em stabilizability} problem, when 
 the existence of at least one stable trajectory is proved by a suitable 
 antinorm~\cite{BS08,  FV12, GLP17, GP13,  GZ15, SDP08}. 
 Since a nontrivial antinorm cannot be defined in the entire space, 
 it is applied to systems with an invariant cone, in particular, for positive systems. 
 Dual norms and antinorms naturally appear when considering the dual systems defined by transposed matrices~\cite{M23, PW08, P22}. For applications of antinorms to convex trigonometry 
 in optimal control see~\cite{L19}.    
 
%\enlargethispage{\baselineskip}
 
 \begin{remark}\label{r.50}
{\em Since a positive homogeneous function 
cannot be concave  on the whole space~$\re^d$,  
antinorms  are usually defined on convex cones. 
Some works consider collections of antinorms on partitions 
of~$\re^d$ to convex cones.   An equivalent definition uses Minkowski functionals on star sets~\cite{MR12, R15}.}
\end{remark}

\bigskip 

\begin{center}

\large{\textbf{2. Duality}}	
\end{center}
\bigskip

\begin{center}

\textbf{2.1. Dual antinorms and polar conic bodies }
\end{center}
\bigskip

\begin{defi}\label{d.40}
For an antinorm~$f$ on~$\re^d_+$, its dual antinorm is 
$f^{*}(\by)\, = \, \inf\limits_{f(\bx) \ge 1}(\by, \bx)$. 
\end{defi}
We also write~$f^{*}(\by)\, = \, \inf_{\bx \in \re^d_+}\frac{(\by, \bx)}{f(\bx)}$ 
assuming that the infimum is taken over the points such that~$f(\bx)\ne 0$. 
\begin{defi}\label{d.50}
For a conic body~$G\subset \re^d_+$, its polar is~$\, G^*\, = \, 
\bigl\{\by\in \re^d_+: \ \inf\limits_{\bx \in G}(\by, \bx) \, \ge \, 1\bigr\}$.  
\end{defi}
It follows easily that if an antinorm~$f$ has a unit ball~$G$, then 
the dual antinorm~$f^*$ has the unit ball~$G^*$. The main properties of the duality are 
very similar to those from the classical convex analysis: 
\medskip 

 $f(\bx)\, f^*(\bx) \, \le \, (\bx, \by)$ (an analogue of Young's inequality); 
\smallskip 

  $f^{**} = f; \quad G^{**} = G\ $ (analogues of the Fenchel-Moreau theorem and of the 
  bipolar  
  
  theorem); 
\smallskip 

  $(G_1 \cap G_2)^* \, = \, {\rm co}_s \{G_1^*, G_2^*\}; \qquad 
(G_1 \cup G_2)^* \, = \, ({\rm co}_s \{G_1, G_2\})^* \, = \, G_1^* \cap G_2^*$
 \medskip 
 
\noindent  The first and the third properties follow directly from the definition, 
 the proof of the second one (on the bipolar and double dual), see~\cite{GZ20}. 
 
 Formally, one can define the polar to an arbitrary subset~$M\subset \re^d_+$
 by the same formula: $M^* = \{\by \in \re^d_+: \ \inf_{\bx \in M}(\bx, \by) \ge 1\}$.
 If the closure of~$M$ does not contain the origin, then~$M^*$ is nonempty. 
 Denote by~$G$ the closure of~${\rm co}_+(M)$. Then~$G$ is a conic body and 
 $M^* = G^*$. Thus, the polar to an arbitrary set is a polar to the closure 
 of its positive convex hull (which is conic body). In particular, for every set~$M$ separated from the 
 origin, we have~$M^{**} = \overline{{\rm co}_+(M)}$.  
 
 \begin{figure}[h!]
\centering
 	{\includegraphics[scale = 0.2]{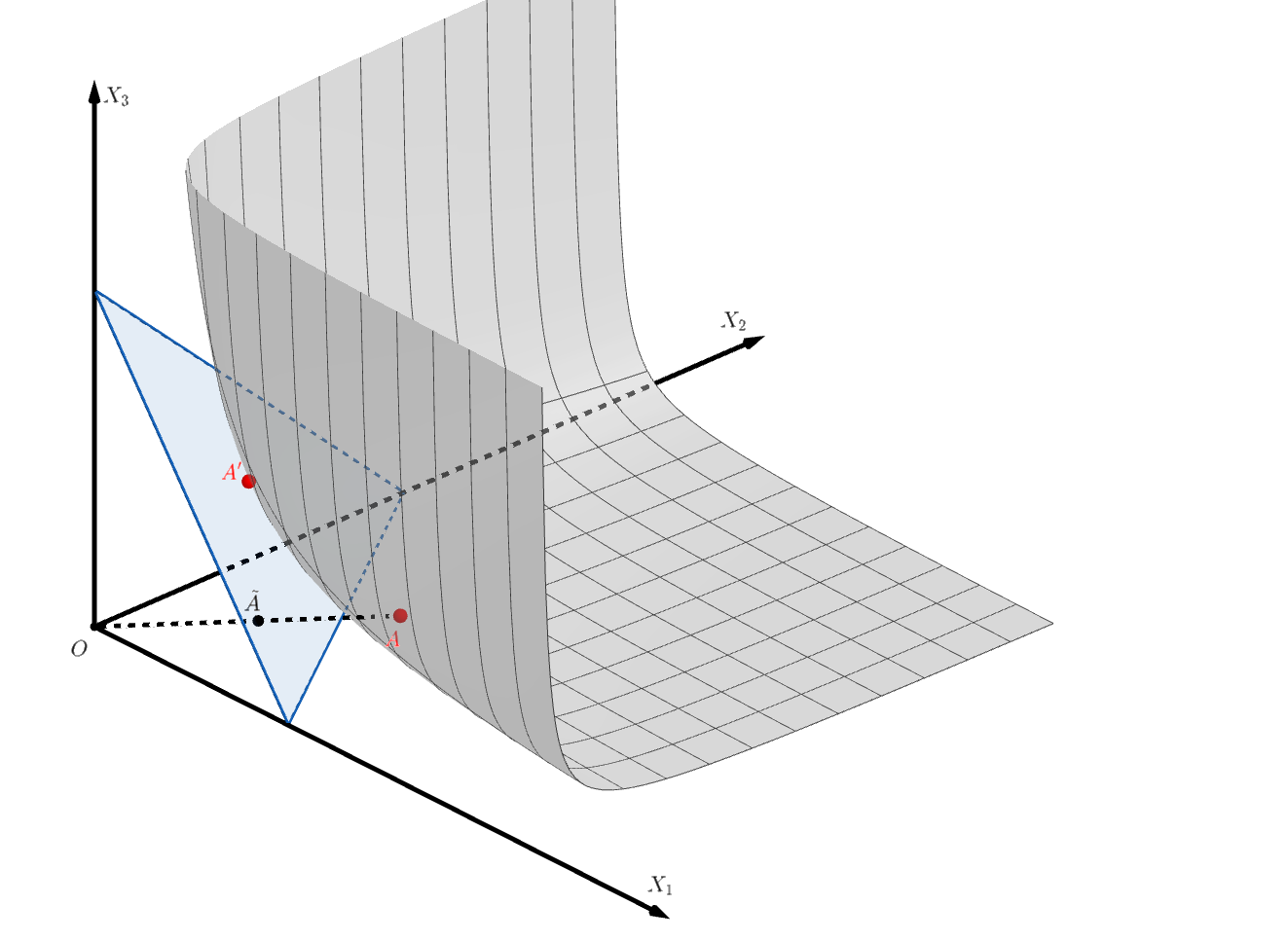}}
 	\caption{The boundary of a conic body in~$\re^3_+$ and the polar to a point~$A$}
 	\label{fig.Smooth-antinorm}
 \end{figure}
 
 \begin{ex}\label{ex.5}
{\em The polar to a point~$\ba \in \re^d_+, \ba \ne 0$, 
 is the same as the polar to the conic polyhedron~${\rm co}_+(\bx) \, = \, 
 \{\bx \in \re^d_+: \, x_i \ge a_i, \, i = 1, \ldots , d\}$ and is 
 equal to the conic polyherdon~$\ba^* \, = \, \{\by \in \re^d_+: \, (\ba, \by) \ge 1\}$. 
 This polyhedron is obtained by cutting off the corner of the positive orthant~$\re^d_+$
  by the hyperplane~$H$ orthogonal 
 to the vector~$\ba$ and passing through 
 the point~$\frac{\ba}{|\ba|^2}$ (inverse to the point~$\ba$). 
 Sometimes we consider this hyperplane is a polar to~$\ba$. 
 
 The polar to a conic polyhedron~$Q$ with vertices~$\ba_1, \ldots \, \ba_n$ 
 (in which case $Q= {\rm co}_+ \{\ba_1, \ldots \, \ba_n\}$)
 is the conic polyhedron~$Q^*$  defined by the system of inequalities~$(\ba_i, \bx) \ge 1, \, 
 i=1, \ldots , n$. }
 \end{ex}
 
 We see that the duality 
theory is very similar to what we have in the classical case, with two exceptions. 
 First, the map~$f\mapsto f^*$ can be discontinuous~\cite{P22}. 
  The second one is the  self-duality issue addressed in the next section.

\bigskip 
\newpage

\begin{center}

\textbf{2.2. Self-duality}	
\end{center}
\bigskip 

An antinorm~$f$  on~$\re^d_+$ is {\em self-dual} if $f^* = f$; a conic body~$G$
is {\em autopolar} if~$G^* = G$. The self-duality of an antinorm is 
equivalent to the autopolarity of its unit ball~$G$. This means that
for an arbitrary point~$A \in \partial \,G$, its polar is a hyperplane of
support for~$G$, and vice versa (Fig.~\ref{fig.Smooth-antinorm}). 

Let us remember that in 
 the convex analysis the self-duality  is very rare. To avoid a confusion we denote the (classical) polar of a set~$Q\subset \re^d$ 
by~$\hat Q$. Thus,~$\, \hat Q\, = \, 
 \{\by \in \re^d: \ \sup\limits_{\bx \in Q} (\by, \bx) \le 1\}$. 
The following fact is well-known, we include its proof for convenience of the reader. 

\smallskip 

\noindent \textbf{Fact.} {\em The Euclidean unit 
ball~$B=\{\bx \in \re^d: |\bx|\le 1\}$ 
is a unique subset of~$\re^d$ such that~$\hat B \, =\, B$}.
\smallskip 

Indeed, if a set~$Q$ coincides with~$\hat Q$, then~$Q \subset B$. 
Otherwise,~$Q$ contains a point~$\by$ such that~$(\by, \by) > 1$ and hence, 
$\by \notin \hat Q$. 
On the other hand, the inclusion~$Q \subset B$ implies~$\hat B \subset \hat Q$ and so, $B \subset Q$, 
because~$Q$ and~$B$ are both autopolar. Hence, $Q=B$. 

\smallskip 

Thus, apart from the unit Euclidean ball, there are no 
autopolar sets in~$\re^d$. 
One, therefore, could expect that there are no autopolar conic 
bodes in~$\re^d_+$ either. However, the situation is quite opposite: there are a lot of autopolar conic bodies, including polyhedra. 
This phenomenon was first observed for~$d=2$ in~\cite{P22}. 
For general~$d \ge 3$, the orthogonal 
extensions (preimages of projections) of a planar autopolar body is also autopolar. 
 Are there other examples?  
The  following family of smooth autopolar conic bodies  (with smooth strictly convex boundary)
was found in~\cite{P22}. We describe it by the corresponding antinorms: 
  \smallskip 
 
\begin{ex}\label{ex.50}{\em 
For an arbitrary collection of non-negative numbers $\{p_i\}_{i=1}^d$ 
such that   $\sum_{i=1}^d p_i = 1$, the function  
\begin{equation}\label{eq.prod}
f(\bx) \ = \  \prod_{i=1}^d\,  \left( \frac{x_i}{\sqrt{p_i}}\right)^{\,p_i} \ , \qquad \bx \in \re^d_+
\end{equation}
(for~$p_i=0$ we set~$\bigl( \frac{x_i}{\sqrt{p_i}}\bigr)^{\,p_i} = 0$) is a 
continuous self-dual antinorm.   
For example, if $p_1=p_2 = 1/2$, then we obtain a self-dual antinorm 
$f(x_1, x_2) = \sqrt{2x_1x_2}$ in~$\re^2_+$. The corresponding autopolar 
ball~$G$ is bounded by the hyperbola. 

The function~(\ref{eq.prod}) can be rewritten in a simpler form 
as~$f(\bx) = Cx_1^{p_1}\cdots x_d^{p_d}$, where~$C$ is a constant 
depending on~$p_i$. There was a mistake in computing this constant in~\cite{P22}, 
now we correct it. 
In Proposition~\ref{p.30} we show that the self-duality of~(\ref{eq.prod}) 
follows from the general result of Theorem~\ref{th.20} along with infinitely many 
other smooth  self-dual antinorms.  
}
\end{ex}

This was the only known example of nontrivial smooth self-affine antinorm. 
Theorem~\ref{th.20} from the next section  provides a rich variety of 
 other examples. 
 
The problem of existence of nontrivial  autopolar conic polyhedra in dimensions~$d\ge 3$
 was left in~\cite{P22} and was positively resolved in~\cite{M23}.   
For every~$d\ge 3$ and~$n\ge 1$, there exists an autopolar conic polyhedron  
of $n$ vertices in~$\re^d_+$.  We consider one example for~$d=3, n=1$.

\begin{ex}\label{ex.57}{\em In the positive orthant~$\re^3$ we 
take an arbitrary point~$A_1$ on the axis~$OX_3$ such that~$OA_1 > 1$
(Fig.~\ref{fig.P3-1}). 
\begin{center}
\begin{figure}[h!]
\centering
 	{\includegraphics[scale = 0.7]{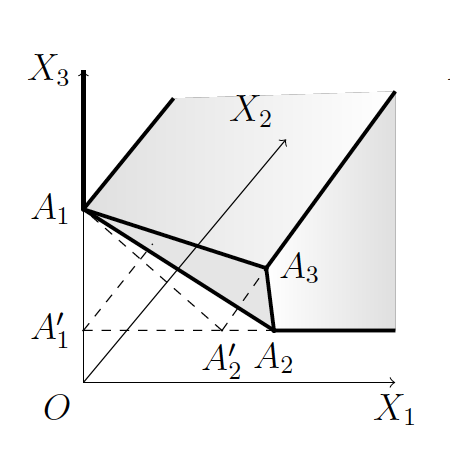}}
 	\caption{The autopolar conic polyhedron~$P$}
 	\label{fig.P3-1}
 \end{figure}
 \end{center}
Denote by~$\alpha_1$ the polar (hyperplane) to~$A_1$. 
It is parallel to~$OX_1X_2$ and passes through the point~$A_1' \in OX_3$ 
such that~$OA_1' = \frac{1}{OA_1}$. Denote by~$\ell_2$ the ray of 
intersection of~$\alpha_1$ and 
the angle~$OX_1X_3$. Take an arbitrary point~$A_2 \in \ell_2$ and denote 
by~$\alpha_2$ its polar (hyperplane).  Let it intersect~$\ell_2$ at a point~$A_2'$. 
Denote by~$\ell_3$ and~$\ell_1$ the rays of intersection of~$\alpha_2$ 
with~$\alpha_1$ and with~$OX_2X_3$ 
respectively. Note that the point~$A_1$ 
is the end of~$\ell_1$. Finally, take a point~$A_3$ on~$\ell_3$ such that 
the three angles~$A_1A_3O, \, A_2A_3O$ and~$A_1A_3A_2$ are all right. 
Such a point exists and unique since the ray~$\ell_3$ is a perpendicular 
to the plane of the triangle~$A_1A_2O$ erected at its ortocenter. Denote by~$\alpha_3$
the hyperplane of the triangle~$A_1A_2A_3$. 
Then the conic polyhedron~$P={\rm co}_+\, \{A_1, A_2, A_3\}$ is autopolar. 
Indeed, it is bounded by three hyperplanes (not counting the coordinate planes)
$\alpha_1, \alpha_2, \alpha_3$ and each plane~$\alpha_i$ is a polar to~$A_i, \, i = 1, 2, 3$. 
See~\cite{M23} for the details. 

The polyhedron~$P$ has three vertices~$A_1, A_2, A_3$, three 
facades~$A_1A_2A_3,\  \ell_2A_2A_3\ell_3, \ \ell_1A_1A_3\ell_3$, and 
three facades on the coordinate subspaces.

}
\end{ex}

No 
autopolar polyhedra different from those presented in~\cite{M23} were unknown.  
Theorem~\ref{th.20} gives an infinite variety 
of autopolar bodies both smooth and polyhedral. Before formulating it we
observe several general properties of autopolar objects.  
\smallskip

\begin{prop}\label{p.40}
If antinorms~$f$ and~$g$ are self-dual and~$f\ge g$, then~$f=g$. 
\end{prop} 
{\tt Proof}. The definition of dual antinorm implies that 
if~${f\ge g}$, then~${f^*\le g^*}$ and therefore~${f\le g}$.

 \smallskip 

{\hfill $\Box$}
\medskip

\begin{prop}\label{p.20}
The distance from every autopolar conic body to the origin 
is equal to one. 
\end{prop} 
{\tt Proof}. Let~$f$ be a self-dual antinorm generated by an
autopolar body~$G$.  Denote by~$A$ 
the closest to the origin point of~$G$ and by~$\ba$ the vector~$OA$. 
The analogue to Young's inequality (see above) 
yields~$(\ba, \ba) \ge f(\ba)f^*(\ba)$. On the other hand,  
$\ba \in \partial G$ and hence, $f(\ba) = 1$, so~$f^*(\ba) =1$. 
Thus,~$(\ba, \ba) \ge 1$. The hyperplane~$H$
passing through~$A$ orthogonal to~$\ba$ is a plane of support for~$G$. 
Therefore, its polar $H^*$, which is the point~$\frac{1}{|\ba|^2}\, A$ lies 
on the unit sphere of the antinorm~$f^* = f$.  
Hence, the length of $\frac{\ba}{|\ba|^2}$ is at least one, therefore~$|\ba| \le 1$ and
thus, $|\ba| =1$. 

 \smallskip 

{\hfill $\Box$}
\medskip

\begin{ex}\label{ex.58}{\em 
For the unit ball of the smooth antinorm~(\ref{eq.prod})
from Example~\ref{ex.50}, the closest point to the 
origin is~$\ba = (\sqrt{p_1}, \ldots , \sqrt{p_d})$. Clearly, $|\ba| = 1$. 

For the conic polyhedron~$P$
from Example~\ref{ex.57}, the closest point to the 
origin is~$\ba = A_3$ and~$|\ba| = 1$. Indeed, since the plane~$A_1A_2A_3$ is a 
polar to~$A_3$, it follows that~$\ba = |OA_3| = 1$. 
Now Proposition~\ref{p.20} implies that~$A_3$ is the closest point.  
 
} 

\end{ex}

Applying Proposition~\ref{p.20} to a unit ball of a self-dual antinorm, we see that
$|\bx| \ge 1$, whenever $f(\bx) \ge 1$. Therefore~$f(\bx) \ge |\bx|$ and 
this inequality is strict unless~$\bx$ is not collinear to~$\ba$. Thus, we obtain 
\begin{cor}\label{c.15}
Every self-dual antinorm does not exceed the Euclidean norm and there is 
a unique point where they both equal to one.
\end{cor}

\bigskip 

\newpage

\begin{center}

\large{\textbf{3. The main results}}	
\end{center}
\bigskip 

We begin with Theorem~\ref{th.20} that asserts that 
every self-dual antinorm can be extended ({\em lifted}),  from a suitable hyperplane to 
the whole cone~$\re^d_+$. This lifting is not unique and  
every~$(d-1)$ dimensional self-dual antinurm can produce infinitely many 
$d$-demensional ones. This gives a large variety of autopolar conic bodies, both smooth 
and polyhedral, built inductively. Then we prove Theorem~\ref{th.25}  that establishes the structure of those antinorms. For the case $d=2$, this method gives all autopolar sets (Theorem~\ref{th.30}) while for~$d=3$ there is an autopolar polyhedron not obtained this way
(Example~\ref{ex.100}). 

We call a hyperplane~$V \subset \re^d$ {\em admissible} if it has the 
form~$V\, = \, \bigl\{\bx \in \re^d: \ x_i = \mu x_j \bigr\}$, where 
$i,j$ are two different indices from~$\{1, \ldots , d\}$ and~$\, \mu \ge 0$ is
arbitrary. Thus, an admissible hyperplane contains $d-2$ coordinate axes 
and intersect~$\re^d_{+}$ by an~$(d-1)$ dimensional orthant~$\re^d_+$, which will be denoted as~$K'$. 
\medskip

\begin{center}

\textbf{3.1. The lifting theorems}	
\end{center}
\bigskip 

Consider an arbitrary admissible hyperplane $V$ that 
splits the orthant $\re^d_+$ into parts~$K_1, K_2$. 
We denote~$K'=  \re^d_+ \cap V$. 

Let~$\varphi$ be an antinorm 
on~$K'$, $G' \subset K'$ be its unit ball, $\Phi$ be its {\em orthogonal extension} to~$\re^d_+$ defined as~$\Phi(\bx) = \varphi(\bx'), \, \bx\in \re^d_+$, where~$\bx'$
is the orthogonal projection of~$\bx$ to~$V$. Thus, $\varphi$ is defined on~$K'$ while 
$\Phi$ is defined on~$\re^d_+$. The unit ball of~$\Phi$ is an intersection of~$\re^d_+$ with a right cylinder based on~$G'$. 

\begin{theorem}\label{th.20}
Let an admissible hyperplane $V$ split~$\re^d_+$ into closed parts~$K_1, K_2$ and let 
$\varphi$ be a self-dual antinorm on~$K'$. Let~$f_1: K_1\to \re_+$ 
be an arbitrary  antinorm 
 which coincides with~$\varphi$ on~$K'$
and does not exceed~$\Phi$ on~$K_1$. Then 
the function~$f(\bx)$ defined on~$\re^d_+$ by the formula 
 \begin{equation}\label{eq.f}
f(\bx) \ = \ 
\left\{
\begin{array}{lcl}
f_1(\bx)\, & , & \bx \in K_1,\\ 
${}$ & ${}$ & ${}$\\
\inf\limits_{\bz \in K_1}\frac{(\bx, \bz)}{f_1(\bz)}\, & , &\bx \notin K_1. 
\end{array}
\right. 
\end{equation}
 is a self-dual antinorm on~$\re^d_+$. 
\end{theorem}
Thus, a self-dual antinorm~$\varphi$ on~$\re^{d-1}_+$ produces 
infinitely many such antinorms on~$\re^d_+$. For this,~$\varphi$ is placed 
on an arbitrary admissible hyperplane~$V$ and  then 
an arbitrary 
antinorm~$f_1$ not exceeding the orthogonal extension of~$\varphi$ can be taken on the part~$K_1$. 
Then we get a self-dual antinorm by formula~(\ref{eq.f}). 
The structure of this antinorm is clarified by Theorem~\ref{th.25} below. 

Let an admissible  hyperplane  $V$ cut the positive orthant~$K=\re^d_+$
 into two closed parts~$K_1, K_2$ and~$K'= K\cap V$. 
Suppose functions~$f_1, f_2$ are defined on~$K_1, K_2$ respectively 
and coincide on~$K'$; then their  {\em concatenation}~$f = f_1\cup f_2$ 
is defined on~$\re^d_+$ as follows:  
\begin{equation}\label{eq.conc}
f(\bx) \ = \ 
\left\{
\begin{array}{l}
f_1(\bx)\, , \, \bx \in K_1,\\ 
f_2(\bx)\, , \, \bx \in K_2\, . 
\end{array}
\right. 
\end{equation}

\begin{theorem}\label{th.25}
Under the assumptions of~Theorem~\ref{th.20}, define the function
\begin{equation}\label{eq.extension}
f_2(\by) \ = \ \inf_{\bx \in K_1}\frac{(\by, \bx)}{f_1(\bx)}\, , \qquad \by \in K_2.
\end{equation}
 Then

1) $\ f_2|_{V}\, = \, \varphi\, $ and  $\ f_2 \le \Phi$ on~$K_2$; 

2) $\ f_1(\bx) \, = \, \inf_{\by \in K_2}\frac{(\bx, \by)}{f_2(\by)}$; 

\end{theorem}
Thus, $f_1$ and~$f_2$ are obtained from each other by the same formula~(\ref{eq.extension}). Theorem~\ref{th.20} asserts that the concatenation~$f = f_1 \cup f_2$ is a self-dual antinorm on~$\re^d_+$. 

\begin{remark}\label{r.65}
{\em The proofs of Theorems~\ref{th.20} and \ref{th.25} 
are given in Section~5. We shall see that
 the admissibility of the plane~$V$ is used only once: 
 the projection of an arbitrary point~$\bx \in \re^d_+$ to~$V$ lies in~$K'$.  
 We are not aware whether it is possible to 
 generalize Theorem~\ref{th.20}, maybe in somewhat different form,
  to an arbitrary hyperplane~$V$ intersecting 
 the interior of~$\re^d_+$? Example from Section~6 (Fig.~\ref{fig.P5-trans}) shows that 
 the answer can be affirmative.  

}
\end{remark}

\begin{center}

\textbf{3.2. Geometrical formulations}	
\end{center}
\bigskip

Reformulation of Theorems~\ref{th.20} and~\ref{th.25} in terms of conic bodies is straightforward. For an admissible  hyperplane~$V$ that splits~$\re^d_+$ to the parts~$K_1, K_2$, (for the sake of simplicity, we assume that they both possess a nonempty interior), 
denote~$K'=\re^d_+ \cap V$. 

\begin{cor}\label{c.10}
Let~$V$ be an  admissible hyperplane, 
  $G'$ be a $(d-1)$-dimensional autopolar conic body in~$K'$. Then for an arbitrary
conic body~$G_1 \subset K_1$ whose intersection with~$V$ 
and projection to~$V$ coincide with~$G'$, the set~$G \, = \, G_1\cup G_2$, where 
 \begin{equation}\label{eq.G}
G_2 \ = \ \{\by \in K_2: \inf_{\bx \in G_1} (\by, \bx) \ge 1\}
\end{equation}
is an  autopolar conic body in~$G$.    
\end{cor}
Below we use the notation from Corollary~\ref{c.10}: 
the conic body~$G_2 \subset K_2$ is obtained from~$G_1$ by formula~(\ref{eq.G})
and~$G = G_1\cup G_2$. 
Let us recall that an orthogonal extension of a set~$M\subset V$ is the 
preimage of~$M$ of the orthogonal projection~$\re^d \to V$. 
\begin{cor}\label{c.20}
If  $G'$ is autopolar in~$K'$ and the orthogonal extension of every $(d-2)$-dimensional plane of support 
of~$G'$ is a $(d-1)$-dimensional plane of support for~$G_1$, then~$G$
is an autopolar conic body.   
\end{cor}

If all those planes of support are tangents to~$S_1 = \partial G_1$ and~$S_1$ is smooth and strictly convex,
 then so is~$S$.   This way one can produce smooth and strictly convex autopolar surfaces. 
Now turn to conic polyhedra. 

\begin{cor}\label{c.30}
If~$G_1 \subset K_1$ is an arbitrary  conic polyhedron, its  facade~$G' = G_1\cap V$   
is autopolar in~$K'$ and all dihedral angles  adjacent 
to this facade do not exceed~$90^{\circ}$, then $G$ is an autopolar polyhedron.
\end{cor}

Let~$G'$ be given by a system of~$n$ linear inequalities~$(\ba_i, \bx) \ge 1\, \, \bx \in K'$, 
with all~$\ba_i \in K'$. We add arbitrary linear inequalities~$(\ba_j, \bx) \ge 1, \ j\in n+1, \ldots , n+m; \ \bx \in \re^d_+$ such that~$\ba_{j} \in K_2$
and the hyperplanes~$(\ba_j, \bx) = 1$ do not intersect the relative interior of~$K'$. 
Then the conic polyhedron~$G_1 \, = \, \{\bx \in K_1: \ (\ba_i, \bx) \ge 1, \ 
i=1, \ldots , n+m \}$ satisfies the assumptions of Corollary~\ref{c.30}. 
Conversely, every conic polyhedron~$G_1$ from~Corollary~\ref{c.30} is obtained this way. 
Each of them produces an autopolar polyhedron in~$\re^d_+$. 
\bigskip

\begin{center}

\textbf{3.3. Self-duality in low dimensions}	
\end{center}
\bigskip 

Theorem~\ref{th.20} provides a method of inductive construction 
of autopolar sets starting with the dimension one. 
\medskip

\noindent $\mathbf{d=1}$. This is a trivial case: 
every antinorm on~$K=\re_+$ is  a linear function$f(x) = kx,  \ k> 0$
(we do not use the bold letters since we deal with scalars). 
Then~$f^*(y) = \frac{y}{k}$ and hence, the only self-dual antinorm is~$f(x)=x$. 
\medskip 

\noindent $\mathbf{d=2}$. Theorem~\ref{th.20} produces infinitely many 
self-dual antinorms in~$\re^2_+$. Moreover, it actually gives 
their complete classification.  In this case~$V$ is a line~$\{\, t\, \ba\, : \ t\in \re\}$ 
with~$\ba  \in \re^2_+, |\ba| = 1$, and, 
respectively,~$K'$ is a ray~$\{\,t\, \ba\, : \ t \ge 0\}$. 

Theorem~\ref{th.20} suggests the following
\smallskip

\noindent \textbf{Algorithm~1 of constructing self-dual antinorms in~$\re^2_+$}. 
We take an arbitrary suitable unit vector~$\ba \in \re^2_+$ and set
$\varphi(t\ba) = t, \ t\ge 0$. This is a self-dual antinorm on~$K'$. 
Its orthogonal extension to~$\re^2_+\, $ is~$\, \Phi(\bx) = (\ba, \bx)$.  Then choose an arbitrary antinorm~$f_1$ 
on~$K_1$ (the angle between~$V$ and~$OX_1$ axis) such 
that~$f_1(\ba) = 1$ and~$f_1(\bx) \le (\bx, \ba), \, \bx \in K_1$.
Geometrically this means that the unit ball~$G_1 \subset K_1$
is separated from the origin by the perpendicular to~$V$ erected 
at the point~$\ba$. Equivalently, the tangent to~$S_1$ at the point~$\ba$
 makes an angle of at most~$90^{\circ}$ with the ray~$K'$. 
Then~$f_2$ is defined in~$K_2$ by formula~(\ref{eq.extension})
and ~$f=f_1\cup f_2$ is a self-dual antinorm on~$\re^2_+$. 

\smallskip 

\begin{theorem}\label{th.30}
Every self-dual antinorm on~$\re^2_+$ is obtained by Algorithm~1. 
\end{theorem} 
{\tt Proof}. Let~$f$ be an arbitrary  self-dual  antinorm~on~$\re^2_+$. 
Denote by~$A$
the closest to the origin~$O$ point of the conic sphere~$S = \partial G$, 
$\ba$ is the vector~$OA$. The line~$\ell$ passing through~$A$ 
and orthogonal to~$\ba$ separates~$S$ from~$O$, hence the values of~$f$ on~$\ell$
do not exceed one. By Proposition~\ref{p.20}, $|\ba| = 1$, hence, for all~$\bx \in \ell$, 
$f(\bx) \le 1 \, = \, (\ba, \bx)$. Therefore, by homogeneity, 
$f(\bx) \le (\ba, \bx)$ for all~$\bx \in \re^2_{+}$. 
 The ray~$K'$ splits~$\re^2_+$ into two angles
$K_1$ and $K_2$ and the function~$f_1 = f|_{K_1}$
satisfies all the assumptions of~Theorem~\ref{th.20}. 
Hence, the function~$\tilde f = f_1 \cup f_2$, where 
$f_2(\by) = \inf\limits_{\bx \in K_1} \frac{(\by, \bx)}{f(\bx)}$
is self-dual. On the other hand, for every~$\by \in K_2$, we have 
$$
f(\by) \ = \ f^*(\by) \ = \ 
\inf\limits_{\bx \in \re^d_2} \frac{(\by, \bx)}{f(\bx)} \ \le \ 
\inf\limits_{\bx \in K_1} \frac{(\by, \bx)}{f(\bx)} \ = \ \tilde f(\by)\, .  
$$
Thus, $\, f\, \le \tilde f$. Since those functions are both self-dual, we get
$f=\tilde f$ (Proposition~\ref{p.40}).

 \smallskip 

{\hfill $\Box$}
\medskip

 \begin{cor}\label{c.80}
Every autopolar conic polygon in~$\re^2_+$ has the form~$G=G_1\cup G_2$,
where~$K'$ is an arbitrary ray splitting~$\re^2_+$ to two angles,~$G_1$ 
is an arbitrary conic polygon in~$K_1$
whose angle at the  vertex~$\ba$ is at most~$90^{\circ}$
and $G_2$ is a polar transform of~$G_1$ obtained by~(\ref{eq.G}). 
 \end{cor}
 
 \bigskip

\begin{center}

\large{\textbf{4. Examples and special cases}}	
\end{center}
\bigskip 

We consider several examples of applications of Theorems~\ref{th.20} and~\ref{th.25} 
for constructing autopolar sets in general dimensions~$d\ge 2$. 

\medskip

\begin{ex}\label{ex.55}
\textbf(The case of right cylinder). {\em In Theorem~\ref{th.20} one can always take~$f_1 = \Phi$, this is the maximal possible~$f_1$. In this case~$G_1$ is a right cylinder with the base~$G'$. Let~$V$ separates the axes~$OX_{d-1}$ and~$OX_d$ and~$K_1$ contains the first of them. 
The function~$f_2$ and the set~$G_2$ is found by Theorem~\ref{th.25}. Its 
unit ball~$G_2$ is an oblique  cylinder with base~$G'$ and with the 
element parallel to the~$OX_d$ axis. The conic body~$G$ composed of those two cylinders is autopolar.  
}
\end{ex}

\begin{ex}\label{ex.60}
\textbf(The orthogonal extension). {\em In the limit case the hyperplane~$V$ 
coincides with one of the coordinate planes, say, the plane~$OX_1\ldots OX_{d-1}$. 
In this case $f_1$ must be equal to~$\Phi$. 
Indeed, if~$f_1(\bx) < \Phi(\bx)$ for some~$\bx$, 
then~$f_1(\bx' + x_d\be_d) < \varphi(\bx')$, where~$\be_d$ is the $d$th basis vector
and~$\bx'$ is the projection of~$\bx$ to~$V$. The  concavity of the function~$f_1$
implies that it becomes negative for large~$x_d$, which is impossible. 

Thus, if~$V$ coincides with a coordinate plane, then~$f_1 = \Phi|_{K_1}$ and the only feasible 
conic body~$G_1$ is an orthogonal extension of~$G'$. Hence,
the autopolar set~$G$ produced by Theorem~\ref{th.20} is an orthogonal extension of~$G'$.  
}
\end{ex}

\begin{ex}\label{ex.65}
\textbf(Autopolar polygons) {\em  In view of Theorem~\ref{th.30}, every autopolar 
conic polygon is constructed as follows. We take a point~$A_0 \in \re^2_+$
such that~$OA_0 = 1$ and choose an arbitrary conic polygon~$P = A_0A_1\ldots A_n$ inside the 
angle $A_0OX_1$ such that the angle of~$P$ at the vertex~$A_0$ is at most~$90^{\circ}$. 
The polygon~$P$ has two infinite sides: one starting at the vertex~$A_0$ parallel to~$OA$
(denote it by~$a$), the other  starting at~$A_n$ parallel to~$OX_1$. This is~$G_1$. Denote by~$a_i$ the line which is the
polar to~$A_i$. The lines~$a, a_0, \ldots , a_n$ (in this order) form the polygon~$G_2$
in the angle~$AOX_2$. The polygon~$G=G_1\cup G_2$ is autopolar (Fig.~\ref{fig.2D}). 
\begin{figure}[h!]
\centering
 	{\includegraphics[scale = 0.6]{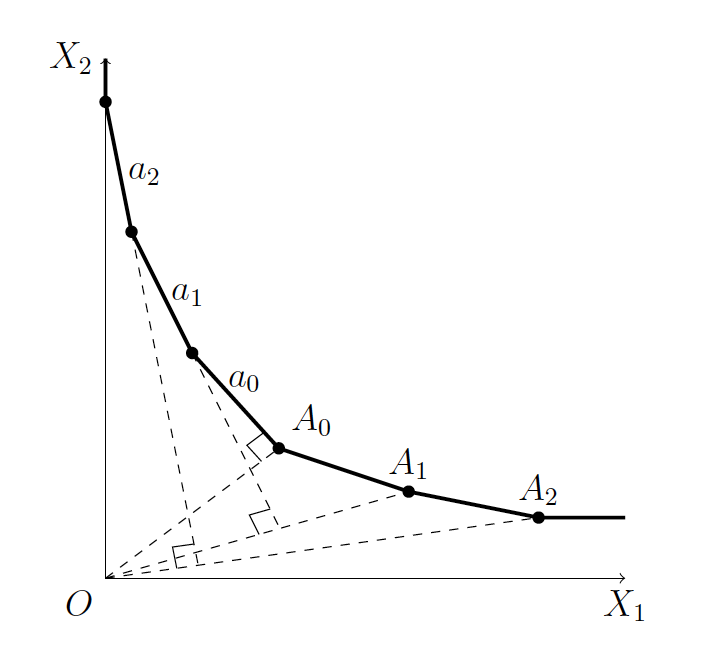}}
 	\caption{An autopolar conic polygon}
 	\label{fig.2D}
 \end{figure}
 
}
\end{ex}

Now turn to the smooth antinorm from~Example~\ref{ex.50} and show that 
it is also constructed by Theorem~\ref{th.20}. 

\begin{prop}\label{p.30}
For every non-negative numbers $p_1, \ldots , p_{d}$ 
such that   $\sum_{i=1}^d p_i = 1$, the self-dial antinorm~(\ref{eq.prod}) 
is obtained by Theorem~\ref{th.20} from the $(d-1)$-dimensional 
  self-dual antinorm~$\varphi$ corresponding to the 
  numbers~$p_1, \ldots , p_{d-2}, p_{d-1}  + p_{d}$. 
\end{prop} 
{\tt Proof.} 
We assume that all~$p_i$ are nonzero, 
otherwise we reduce the dimension. Let us show that~$f$ 
is obtained by Theorem~\ref{th.20} with the  
hyperplane~$V\, = \, \bigl\{\bx \in \re^d \ : \ 
\frac{x_{d-1}}{\sqrt{p_{d-1}}} \, = \, \frac{x_d}{\sqrt{p_d}} \, \bigr\}$
and~$f_1 = f|_{K_1}$.  
Introduce the  coordinate system in~$V$ as follows: 
all the axes~$OX_i, \, i = 1, \ldots, d-2$ are the same as in~$\re^d$
and the axis~$OX_0 = (OX_{d-1}X_{d})\cap V$. We have  
$$
x_{d-1} \ = \ \frac{\sqrt{p_{d-1}}}{\sqrt{p_{d-1}+ p_d}} \ x_0\, , \qquad 
x_{d} \ = \ \frac{\sqrt{p_{d}}}{\sqrt{p_{d-1}+ p_d}} \ x_0
$$
Denote~$p_0 = p_{d-1}+p_{d}$ and consider  the function~$\varphi(\bx)\, = \, 
\prod\limits_{i=0}^{d-2}\,  \left( \frac{x_i}{\sqrt{p_i}}\right)^{\,p_i} 
$ on~$K'$. By the inductive assumption, 
$\varphi$ is self-dual. Then the gradient 
$$
f'(\bx)\ = \ \left( \frac{p_1}{x_1}\, , \ldots , \frac{p_d}{x_d} \right) \, f(\bx)\, . 
$$
belongs to~$V$, whenever $\bx \in V$. Indeed, 
if~$\frac{x_{d-1}}{x_{d}} \, = \, \frac{\sqrt{p_{d-1}}}{\sqrt{p_d}}$, then 
$f_{x_{d-1}}/f_{x_d} \, = \, 
\frac{p_{d-1}}{x_{d-1}}f(\bx)\, / \,   \frac{p_d}{x_d}f(\bx)\, = \, \frac{\sqrt{p_{d-1}}}{\sqrt{p_d}}$. 
Thus, the tangent plane to~$S$ at every point~$\bx \in V\cap S$, 
is orthogonal to~$S$. In view of Corollary~\ref{c.20}, we have $f\le \Phi$. 

 \smallskip 

{\hfill $\Box$}
\medskip 

For~$d=1$, the antinorm~(\ref{eq.prod}) becomes~$f(x) = x$, which is self-dual. 
Then every $d$-dimensional antinorm of type~(\ref{eq.prod}) can be obtained from this one by $d$ steps of Proposition~\ref{p.30} with preserving the self-duality by Theorem~\ref{th.20}.

\begin{ex}\label{ex.67}{\em  The three-dimensional autopolar conic polyhedron~$P$ 
presented in Example~\ref{ex.57} can be  obtained by the lifting form Theorem~\ref{th.20}. In this case the 
hyperplane~$V$ passes through the axes~$OX_3$ and through the point~$A_3$, see 
Fig.~\ref{fig.P3-proof}. 
\begin{figure}[h!]
\centering
 	{\includegraphics[scale = 0.6]{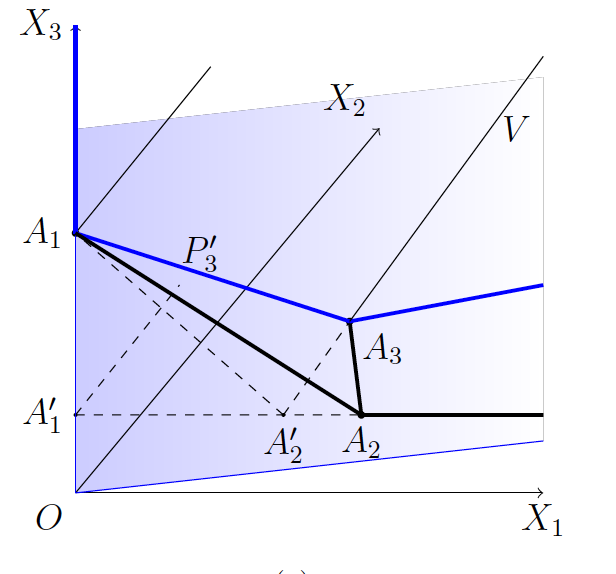}}
 	\caption{The polyhedron~$P$ form Example~\ref{ex.57} is obtained by the lifting}
 	\label{fig.P3-proof}
 \end{figure}
The intersection~$P\cap V$ is a conic polygon~$P'$
bounded by the segment~$A_1A_3$, by the ray~$A_1X_3$, and by the ray~$\ell$ going from~$A_3$
along  the line of intersection of~$V$ with~$\alpha_1$. 
The polyhedron~$P'$ is autopolar. Indeed, the polar to~$A_1$ is the plane~$\alpha_1$, 
hence,~$P'$ is obtained by Algorithm~1 with the partition of~$\re^2_+$ by the line~$OA_3$. 
Then we pass from~$P'$ to~$P$ by Corollary~\ref{c.30}. The assumptions of the 
corollary are satisfied. Indeed,~$V$ is orthogonal to both~$\alpha_1$ and~$\alpha_3$ and hence, 
the dihedral angles adjacent to the edges~$\ell$ and~$A_1A_2$ are right.  
 
}
\end{ex}

\begin{center}

\large{\textbf{5. Proofs of the main theorems}}	
\end{center}
\bigskip 

{\tt Proof of Theorem~\ref{th.25}}.  
For every pair of points~$\bx\in K_1, \, \by \in K'$,  
we have~$(\by , \bx) \, = \, (\by, \bx')$ and
$f_1(\bx) \le \Phi(\bx)   = \varphi(\bx')\, = \, f_1(\bx')$. Therefore, 
$$
\frac{(\by, \bx)}{f_1(\bx)} \ \ge \ \frac{(\by, \bx')}{f_1(\bx')}\, . 
$$ 
Take the infimum of the both parts over~$\bx \in K_1$. 
In the left hand side we obtain~$f_2(\by)$,  the right hand side becomes 
$$
\inf_{\bx\in K_1}\frac{(\by, \bx')}{f_1(\bx')} \ = \ 
\inf_{\bx'\in K'}\frac{(\by, \bx')}{\varphi(\bx')} \ = \ 
\varphi^* (\by)\, . 
$$
Thus, for every~$\by \in K'$, we have~$f_2(\by) \ge  \varphi^*(\by)$. 
On the other hand, an infimum over the set~$K_1$ does not exceed an infimum 
over~$K'$, 
hence~$f_2(\by) \, \le  \, \inf\limits_{\bx\in K'}\frac{(\bx, \by)}{f_1(\bx)} \, =\, 
\, \varphi^*(\by)$. Thus, $f_2(\by) = \varphi^* (\by)$ and therefore, 
$f_2|_{V} = \varphi^* = \varphi$.

Now take arbitrary~$\by \in K_2$. For every~$\bx \in K_1$, 
we have 
$$
(\by, \bx)\, - \, (\by', \bx) \ = \ (\by - \by', \bx) \ = \ 
(\by - \by', \bx') \, + \, (\by - \by', \bx - \bx') \ \le \  0\, , 
$$
because~$(\by - \by', \bx') = 0$ and 
$(\by - \by', \bx - \bx') \, \le \,  0$ since the vectors~$\by - \by'$ and~$\bx - \bx'$ 
have opposite directions. Thus, $(\by, \bx)\, - \, (\by', \bx)\, \le \, 0$. 
Consequently, 
$$
\frac{(\by, \bx)}{f_1(\bx)}  \ \le \  \frac{(\by', \bx)}{f_1(\bx)}  
$$
Taking an infimum of both parts  over~$\bx \in K_1$, 
we conclude that~$f_2(\by) \le f_2(\by') = \varphi(\by')$, which 
proves the inequality~$f_2 \le \Phi$ on~$K_2$. The proof of~1) is completed. 

\smallskip

Take arbitrary~$\by \in K_2$ and show that 
for every~$\bu \in K_2$, there exists~$\bx \in K_1$
for which~$\frac{(\by, \bx)}{f_1(\bx)} \, \le \, \frac{(\by, \bu)}{f_2(\bu)}$. 
This will imply that the infimum of the value~$\frac{(\by, \bz)}{f(\bz)}$ 
over~$\bz \in \re^d_+$ (which is~$f^*(\by)$) is  achieved on~$K_1$, i.e., 
is equal to~$\inf\limits_{\bz \in K_1}\frac{(\by, \bz)}{f_1(\bz)}\, = \, f_2(\by)$. 
This proves the equality  $f(\by) \, = \, f^*(\by)$ for~$\by \in K_2$. 
To establish the existence of a desirable~$\bx$ it suffices to take~$\bx = \bu'$, which means 
~$\frac{(\by, \bu')}{f_1(\bu')} \, \le \, \frac{(\by, \bu)}{f_2(\bu)}$. 
This follows from two inequalities:  a)~$f_1(\bu') = f_2(\bu') \ge f_2(\bu)$ due to item~1);  b)
$(\by, \bu') \le (\by, \bu)$ because~$\bu \in K_2$ and  hence, 
the vectors~$\by$ and $\bu - \bu'$ form a non-obtuse angle, so~$0 \le (\by, \bu') \le (\by, \bu - \bu')$. 

Due to item 1), the function $f_2$ on~$K_2$ possess the same properties as~$f_1$
on~$K_1$. 
On the other hand,~$f^*|_{K_2} = f|_{K_2} = f_2$.   Hence, for arbitrary~$\bx \in K_1$, 
we argue as above and show 
that~$f_1^{**}(\bx) \, = \, \inf_{\bu \in K_2}\frac{(\bx, \bu)}{f^*(\bu)}\, = \, 
\inf_{\bu \in K_2}\frac{(\bx, \bu)}{f_2(\bu)}$. Invoking now the equality
$f^{**} = f$, we conclude the proof of item~2).

\smallskip 

{\hfill $\Box$}
\medskip

{\tt Proof of Theorem~\ref{th.20}}.  In the proof of Theorem~\ref{th.25} we showed that 
$f^*|_{K_2} = f|_{K_2} = f_2$. By assertion 2) of Theorem~\ref{th.25}, 
the function $f_1$ is defined by~$f_2$ in the same way as 
$f_2$ is defined by~$f_1$. Interchanging those functions 
 we get~$f^*|_{K_1} = f|_{K_1} = f_1$. Therefore, $f^* = f$. 
From the definition of the dual function it immediately follows that~$f^*$ is 
an antinorm, hence, so is~$f$.

\smallskip 

{\hfill $\Box$}
\medskip

\begin{center}

\large{\textbf{6. A counterexample and open problems}}	
\end{center}
\bigskip

We are giving an example of a three-dimensional conic polyhedron which cannot be 
 obtained by the lifting procedure from Theorem~\ref{th.20}. This will leave open the 
 problem of complete classification of autopolar sets (Problem~1 below). 
 We begin with auxiliary results. 

For every conic body~$G$, the minimal distance to the origin 
is attained at a unique point~$A \in G$.  
This is a simple consequence of the convexity. 
If~$G$  is autopolar, then~$OG=1$ due to Proposition~\ref{p.20} (Section~2).  
It turns out that if~$G$ is obtained by the lifting from an admissible hyperplane~$V$, then~$V$ passes through~$A$. 
\begin{prop}\label{p.75}
Assume that an autopolar conic body~$G\subset \re^d_+$ is obtained by 
the lifting of a $(d-1)$-dimensional conic body~$G'$ from an admissible 
hyperplane~$V$ (Theorem~\ref{th.20}). Then~$G$ and~$G'$ have the same point closest to the origin
and this point belongs to~$V$. 
\end{prop}
{\tt Proof}. Let~$A$ and~$A'$ be the closest  to  the origin   points of the conic bodies~$G$ and~$G'$ respectively. Since those bodies are both autopolar, it follows 
that~$OA=OA'=1$. On the other hand, $A$ and~$A'$ belong to~$G$ and the uniqueness of the 
closest point yields that~$A=A'$. 

\smallskip 

{\hfill $\Box$}
\medskip

Thus, if an autopolar body~$G$ is obtained by Theorem~\ref{th.20}, then the 
splitting hyperplane~$V$ contains  the closest to the origin point~$A\in G$. 
This reduces the search of this hyperplane for a given autopolar set to $\frac{d(d-1)}{2}$ admissible 
hyperplanes passing through~$A$.

The second property that directly follows from Theorem~\ref{th.25} is the following 
analogue of Corollary~\ref{c.30}.  It gives a necessary condition for a conic polyhedron for being obtained by the lifting. 
\begin{prop}\label{p.85}
If a conic polyhedron~$P$ is obtained by Theorem~\ref{th.20}, then 
the admissible hyperplane~$V$ intersects~$P$ by an autopolar 
$(d-1)$-dimensional polyhedron~$P'$ and divides~$P$ into two polyhedra~$P_1, P_2$ with a common facade~$P'$. This common facade makes only non-obtuse dihedral angles with all adjacent facets of~$P_1$ and~$P_2$.
In particular, every facade of~$P$ divided by~$V$ is orthogonal to~$V$.  
\end{prop}

Now we are ready to present the promised counterexample.   We use the construction from~\cite{M23}. 
To simplify the notation, 
 we sometimes replace the polar by its boundary. For example, 
the polar of a point~$\ba$ is the set~$\{\by \in \re^{d}_+: 
(\by, \ba) = 1\}$. Similarly, we call a polar to a segment~$[\ba_1, \ba_2]$ 
the set~$\{\by \in \re^{d}_+: 
(\by, \ba_i) = 1, \ i=1,2\}$. 
\begin{ex}\label{ex.100}{\em  The following series of autopolar 
conic polytopes~$\{P_n\}_{n \ge 3}$ was presented in~\cite{M23}. 
   For every~$n\ge 3$, the polytope~$P_n$ is constructed in~$\re^3_{+}$ as follows. 
We  take an arbitrary point~$A_1$ on the axis~$OX_3$ such that~$|OA_1| > 1$
 Then take an arbitrary point~$A_2$ on the intersection of the polar to~$A_1$
with the hyperplane~$OX_1X_3$, see Fig.~\ref{fig.P5-1}. 
\begin{figure}[h!]
\centering
 	{\includegraphics[scale = 0.7]{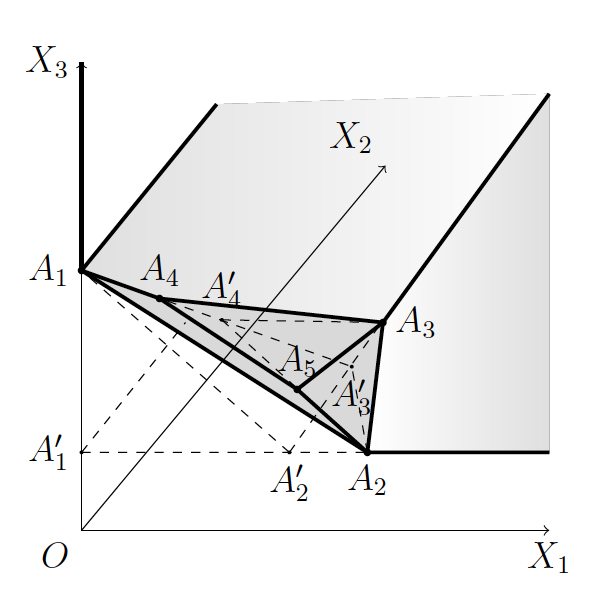}}
 	\caption{Construction of the autopolar conic polyhedron~$P_n$}
 	\label{fig.P5-1}
 \end{figure}
  Then the construction is by induction. 
Suppose the points~$A_1, \ldots , A_k$ are constructed, $k\le n-2$; then~$A_{k+1}$
 is a point on the polar to the segment~$A_{k-1}A_k$ such that~$|OA_{k+1}| > 1$. 
 Finally, on the last iteration, $A_n$ is chosen on the polar to the segment~$A_{n-2}A_{n-1}$
 so that~$|OA_{n}| > 1$. Then~$P_n = {\rm co}_+\, \{A_1, \ldots , A_n\}$
 is an autopolar conic polyhedra, for which~$A_n$ is the point closest to 
 the origin. The proof  is in \cite[theorem~2]{M23}.
 Moreover, the construction is well-defined for every choice of points~$A_1, \ldots , A_n$
 on the corresponding polars under the conditions on the lengths of~$OA_{k}$, see~\cite{M23}. 
 Therefore, we actually construct  a family of autopolar polyhedra and take one of its representatives~$P_n$. Figure~\ref{fig.P5-1} presents the polyhedron~$P_5$, Fig.~\ref{fig.P5-2}
 below shows the same polyhedron without auxiliary lines. 
\begin{center}
\begin{figure}[h!]
\centering
 	{\includegraphics[scale = 0.5]{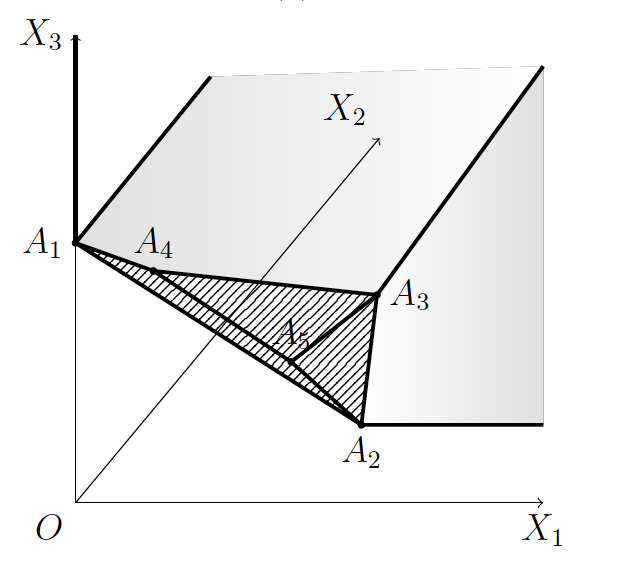}}
 	\caption{The autopolar conic polyhedron~$P_5$}
 	\label{fig.P5-2}
 \end{figure}
 \end{center} 
 The polyhedron~$P$ constructed in Example~\ref{ex.57}   (Fig.~\ref{fig.P3-1}),  
 also belongs to this series,  and actually~$P=P_3$.  As we saw it can be obtained by 
 Theorem~\ref{th.20}. The same is true for~$P_4$. It turns out, however, that 
 $P_5$ is not obtained by that theorem, which gives a desired counterexample. 
}
\end{ex}

\begin{center}
\begin{figure}[h!]
\centering
 	{\includegraphics[scale = 0.4]{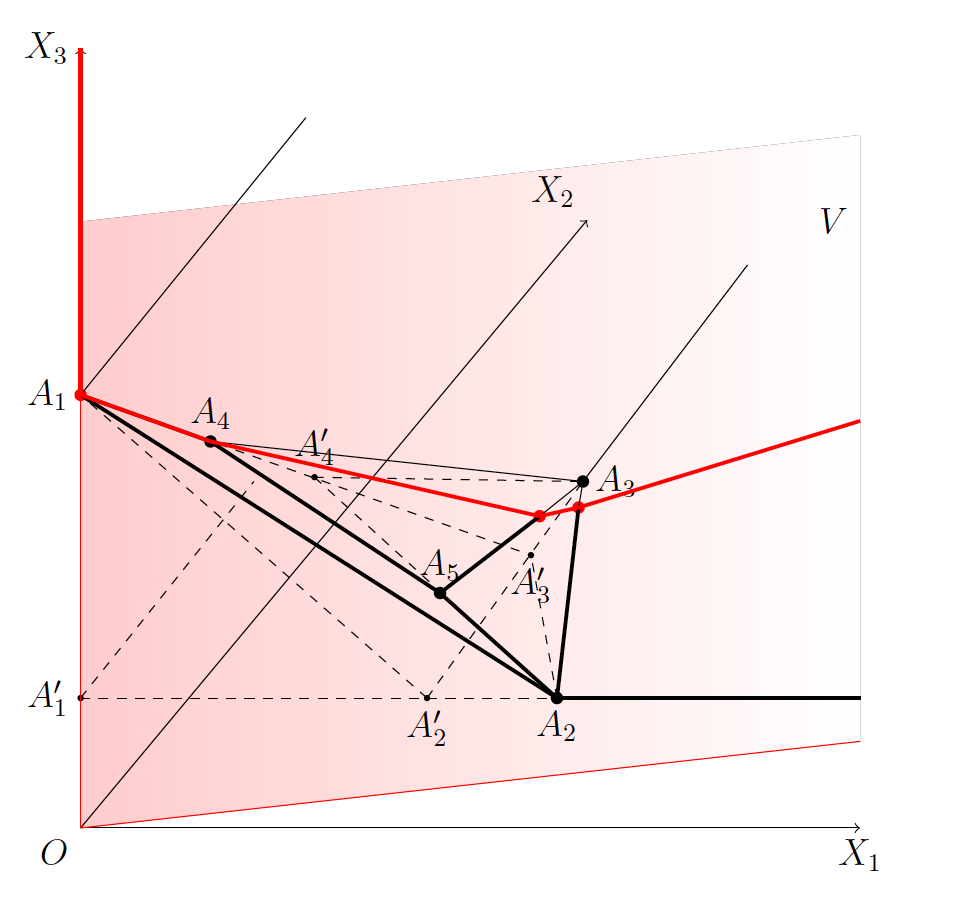}}
 	\caption{The conic polyhedron~$P_5$ cannot be constructed by the lifting}
 	\label{fig.P5-proof}
 \end{figure}
 \end{center}

\begin{center}
\begin{figure}[h!]
\centering
 	{\includegraphics[scale = 0.4]{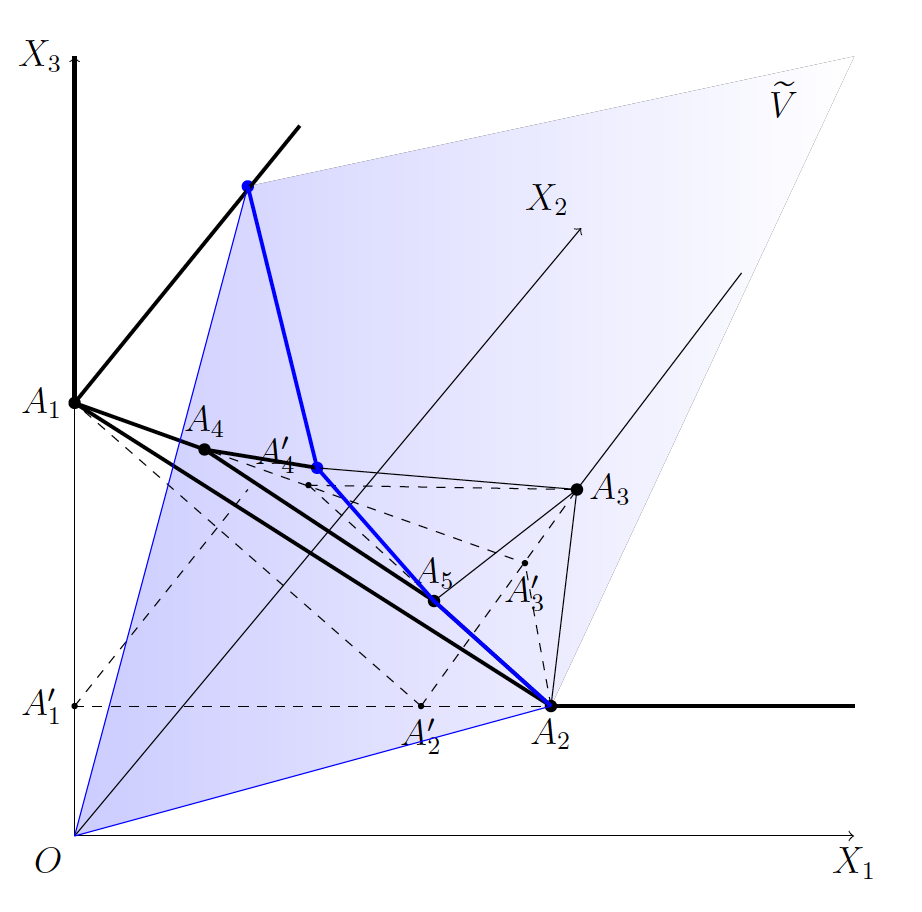}}
 	\caption{The hyperplane~$\tilde V$ orthogonal to all transversal facades of~$P_5$}
 	\label{fig.P5-trans}
 \end{figure}
 \end{center}

\begin{theorem}\label{th.75}
The conic polyhedron~$P_5$ cannot be constructed by the lifting defined in Theorem~\ref{th.20}.
\end{theorem}
{\tt Proof}. If~$P_5$ is constructed by the lifting from an  admissible hyperplane~$V$, then~$V$ passes through the closest to the origin 
point of~$P_5$ (Proposition~\ref{p.75}). This closest point is~$A_5$~\cite{M23}. 
Moreover, due to Proposition~\ref{p.85}, $V$ is orthogonal to all transversal 
facades of~$P_5$. If~$V$ passes through the~$OX_3$ axis, then 
it is transversal (and hence, orthogonal) to the facades~$A_1A_2A_5A_4$ and~$A_2A_3A_5$, 
see Fig.~\ref{fig.P5-proof}. 
Therefore, $V$ is orthogonal to its common edge~$A_2A_5$ and hence, 
$A_2A_5$ is parallel to the coordinate plane~$OX_1X_2$. In this case
the point~$A_5$ lies on the plane~$\alpha_1$, which is the polar to~$A_1$. 
The latter contradicts to the construction of~$P_5$. The cases when~$V$
passes to the other coordinate axes are considered in the same way.

\smallskip 

{\hfill $\Box$}
\medskip

Thus, not all $d$-dimensional autopolar conic polyhedra are obtained by the lifting 
from $(d-1)$-dimensional ones. A natural question arises how to classify all possible autopolar 
sets~? 

\medskip 

\noindent \textbf{Problem~1}. {\em What is the complete classification of 
autopolar conic bodies/polyhedra in~$\re^d_+$?}
\smallskip 

A complete classification can be either a description of all autopolar conic polyhedra 
or an algorithm that can construct all of them. Theorem~\ref{th.30}  solves this problem for~$d=2$,  for 
higher dimensions $d$ it remains open. 
\medskip

A weakened version of Problem~1 is to  generalize 
Theorem~\ref{th.20} to  cover all autopolar  conic bodies. 
For example, does there exist a lifting procedure from a non-admissible hyperplane~$V$, 
i.e., non containing~$(d-2)$ coordinate axes? This can be formulated as follows: 
\medskip 

\noindent \textbf{Problem~2}. {\em Is it true that every autopolar  
conic polyhedron~$P \subset \re^d_+$ is split by a hyperplane~$\tilde V$  
(possibly non-admissible) 
orthogonal to all transversal facades of~$P$ and making non-obtuse angles 
with all facades intersecting it?}
\medskip 

For the polyhedron~$P_5$ such a hyperplane exists: this is the plane~$OA_2A_4$,  see 
Fig.~\ref{fig.P5-trans}. 

In the proof of Theorem~\ref{th.20} we essentially used the fact that 
the orthogonal projection of an arbitrary point~$X\in \re^d_+$ to~$V$ is contained in~$V\subset \re^d_+$. This is true only for admissible hyperplanes~$V$.

 \end{document}